\documentclass{amsart}%
\usepackage{amscd}
\usepackage{amsmath}
\usepackage{amsfonts}
\usepackage{amssymb}
\usepackage{graphicx}%
\providecommand{\U}[1]{\protect\rule{.1in}{.1in}}
\newtheorem{theorem}{Theorem}
\theoremstyle{plain}
\newtheorem{acknowledgement}{Acknowledgement}

\newtheorem{definition}{Definition}

\newtheorem{proposition}{Proposition}

\numberwithin{equation}{section}

\begin{document}
\title[Maximal operators on weighted spaces]{Estimates for Bellman functions related to dyadic-like maximal operators on
weighted spaces }
\author{Antonios D. Melas}
\author{Eleftherios N. Nikolidakis}
\author{Dimitrios Cheliotis}
\address{Department of Mathematics, University of Athens, Panepistimiopolis 15784,
Athens, Greece}
\email{amelas@math.uoa.gr, lefteris@math.uoc.gr, dcheliotis@math.uoa.gr}
\date{September 25, 2015}

\begin{abstract}
We provide some new estimates for Bellman type functions for the dyadic
maximal opeator $\mathbb{R}^{n}$ and of maximal operators on martingales
related to weighted spaces. Using a type of symmetrization principle,
introduced for the dyadic maximal operator in earlier works of the authors we
introduce certain conditions on the weight that imply estimate for the maximal
operator on the corresponding weighted space. Also using a well known estimate
for the maximal operator by a double maximal operators on different measures
related to the weight we give new estimates for the above Bellman type functions.

\begin{acknowledgement}
This research has been co-financed by the European Union and Greek national
funds through the Operational Program "Education and Lifelong Learning" of the
National Strategic Reference Framework (NSRF). ARISTEIA I, MAXBELLMAN 2760,
research number 70/3/11913.

\end{acknowledgement}
\end{abstract}
\maketitle

\section{Introduction}

The dyadic maximal operator on $\mathbb{R}^{n}$ is defined by
\begin{equation}
M{}\,_{d}\phi{}(x)=\sup\left\{  \frac{1}{\left\vert S\right\vert }\int
_{S}\left\vert \phi(u)\right\vert du:x\in S\text{, }S\subseteq\mathbb{R}%
^{n}\text{ is a dyadic cube}\right\}  \label{i1}%
\end{equation}
for every $\phi\in L_{\text{loc}}^{1}(\mathbb{R}^{n})$ where the dyadic cubes
are the cubes formed by the grids $2^{-N}\mathbb{Z}^{n}$ for $N=0,1,2,...$.

As it is well known it satisfies the following weak $L^{p}$ inequality (for
martingales known as Doob's inequality)
\begin{equation}
\left\Vert M_{d}\phi\right\Vert _{p}\leq\dfrac{p}{p-1}\left\Vert
\phi\right\Vert _{p}\label{i3}%
\end{equation}
for every $p>1$ and every $\phi\in L^{p}(\mathbb{R}^{n})$ which is best
possible (see \cite{Burk1}, \cite{Burk2} for the general martingales and
\cite{Wang} for dyadic ones).

An approach for studying more in depth the behavior of this maximal operator
is the introduction of the so called Bellman functions (see \cite{Naz})
related to them which reflect certain deeper properties of them by localizing.
Such functions related to the $L^{p}$ inequality (\ref{i3}) have been
precisely evaluated in\ \cite{Mel1}. Actually defining for any $p>1$
\begin{equation}
\mathcal{B}_{p}(F,f)=\sup\left\{  \dfrac{1}{\left\vert Q\right\vert }\int
_{Q}(M_{d}^{\prime}\phi)^{p}:\operatorname{Av}_{Q}(\phi^{p}%
)=F,\operatorname{Av}_{Q}(\phi)=f\right\}  \label{i4}%
\end{equation}
where $Q$ is a fixed dyadic cube, $R$ runs over all dyadic cubes containing
$Q$, $\phi$ is nonnegative in $L^{p}(Q)$ and the variables $F,f$ satisfy
$0\leq f,f^{p}\leq F$ which is independent of the choice of $Q$ (so we may
take $Q=[0,1]^{n}$) and where the localized maximal operator $M_{d}^{\prime
}\phi$ is defined as in  (\ref{i1}) with the dyadic cubes $S$ being assumed to
be contained in the ambient dyadic cube $Q$. It has been shown in \cite{Mel1}
that
\begin{equation}
\mathcal{B}_{p}(F,f)=F\omega_{p}\left(  \dfrac{f^{p}}{F}\right)
^{p}\label{i5}%
\end{equation}
where $\omega_{p}:$ $[0,1]\rightarrow\lbrack1,\frac{p}{p-1}]$ is the inverse
function of $H_{p}(z)=-(p-1)z^{p}+pz^{p-1}$. Actually (see \cite{Mel1}) the
more general approach of defining Bellman functions with respect to the
maximal operator on a nonatomic probability space $(X,\mu)$ equipped with a
tree $\mathcal{T}$ (see Section 2) can be taken and the corresponding Bellman
function is always the same. The fact that the range of $\omega_{p}$ is
$[1,\frac{p}{p-1}]$ shows in a sense the extend that the constant in Doob's
inequality can be approached only by functions whose integral is very small
compared to its $p$-norm. For example for $p=2$ we get the following sharp
improvement of Doob's inequality
\begin{equation}
\left\Vert M_{\mathcal{T}}\phi\right\Vert _{2}\leq\left\Vert \phi\right\Vert
_{2}+(\left\Vert \phi\right\Vert _{2}^{2}-\left\Vert \phi\right\Vert _{1}%
^{2})^{1/2}<2\left\Vert \phi\right\Vert _{2}\label{i12}%
\end{equation}
which aside from the $L^{2}$ norm of $\phi$ involves also in a sharp way the
\textit{variance }of $\phi$.

Here we will be concerned with the behavior of these maximal operators on
weighted spaces. As it is well known for any positive locally integrable
function $w$ on $Q$ the estimate%
\begin{equation}
\int_{Q}(M_{d}^{\prime}\phi)^{p}w\leq C\int_{Q}\phi^{p}w\label{a}%
\end{equation}
holds for all $\phi$ if and only if $w$ is a dyadic \thinspace$A_{p}$ weight
in the sense that
\[
\sup\{\left\vert I\right\vert ^{-p}(\int_{I}w)(\int_{I}w^{-\frac{1}{p-1}%
})^{p-1}:I\text{ dyadic subcube of }Q\}=[w]_{p}<+\infty\text{.}%
\]
Also it is known that the best possible $C$ is of the order of $[w]_{p}%
^{p/(p-1)}$, the exponent being best possible. Related to this one may define
the following Bellman function given a weight $w$%
\begin{equation}
\mathcal{B}_{p,w}(F,f)=\sup{\huge \{}\dfrac{1}{\left\vert Q\right\vert }%
\int_{Q}(M_{d}^{\prime}\phi)^{p}w:\operatorname{Av}_{Q}(\phi^{p}%
w)=F,\operatorname{Av}_{Q}(\phi)=f{\huge \}}\label{i6a}%
\end{equation}
which is finite only if $w$ is in $A_{p}$ and seek estimates for this in order
to improve the above estimate (\ref{a}). One may add more variables to the
above Bellman function as the integrals of $w$ and of $w^{-1/(p-1)}$ over $Q$
but we will not treat those here. The estimates here will be proved in the
general setting of tree like families on probability spaces and its related
maximal operator, as will be described in the next section.

\bigskip`We will derive two types of estimates related to the above problems.
In the first we will use a related condition on some symmetrization of the
weight to find the exact form of a related to weights Bellman function and
this is done in section 2. Then in section 3 we obtain certain new estimates
for the above Bellman function related to $A_{p}$ with respect to a tree, and
to the corresponding maximal operator, by using an estimate of the maximal
operator via two applications of maximal operators on the same tree but with
different measures, and this is described in section 3.

There are several other problems in Harmonic Analysis where Bellman functions
naturally arise. Such problems (including the dyadic Carleson imbedding and
weighted inequalities) are described in \cite{Naz2} (see also \cite{Naz},
\cite{Naz1}) and also connections to Stochastic Optimal Control are provided,
from which it follows that the corresponding Bellman functions satisfy certain
nonlinear second order PDE.

The exact computation of a Bellman function is a difficult task which is
connected with the deeper structure of the corresponding Harmonic Analysis
problem. Thus far several Bellman functions have been computed (see
\cite{Burk1}, \cite{Burk2}, \cite{Mel1}, \cite{Sla}, \cite{Sla1}, \cite{Vas},
\cite{Vas1}, \cite{Vas2}). L.Slavin and A.Stokolos \cite{SlSt} linked the
Bellman function computation to solving certain PDE's of the Monge Ampere
type, and in this way they obtained an alternative proof of the Bellman
functions relate to the dyadic maximal operator in \cite{Mel1}. Also in
\cite{Vas2} using the Monge-Ampere equation approach a more general Bellman
function than the one related to the dyadic Carleson imbedding Theorem has be
precisely evaluated thus generalizing the corresponding result in \cite{Mel1}.

\bigskip

\section{Trees, maximal operators and symmetrization}

As in \cite{Mel1} we let $(X,\mu)$ be a nonatomic probability space (i.e.
$\mu(X)=1$). Two measurable subsets $A$, $B$ of $X$ will be called almost
disjoint if $\mu(A\cap B)=0$. Then we give the following.

\begin{definition}
A set $\mathcal{T}$ of measurable subsets of $X$ will be called a tree if the
following conditions are satisfied:

(i) $X\in\mathcal{T}$ \ and for every $I\in\mathcal{T}$ \ we have $\mu(I)>0$.

(ii) For every $I\in\mathcal{T}$ \ there corresponds a finite subset
$\mathcal{C}(I)\subseteq\mathcal{T}$ \ containing at least two elements such that:

\qquad(a) the elements of $\mathcal{C}(I)$ are pairwise almost disjoint
subsets of $I$,

\qquad(b) $I=\bigcup\mathcal{C}(I)$.

(iii) $\mathcal{T}=\bigcup_{m\geq0}\mathcal{T}_{(m)}$ where $\mathcal{T}%
_{(0)}=\{X\}$ and $\mathcal{T}_{(m+1)}=\bigcup_{I\in\mathcal{T}_{(m)}%
}\mathcal{C}(I)$.

(iv) We have $\lim\limits_{m\rightarrow\infty}\sup\limits_{I\in\mathcal{T}%
_{(m)}}\mu(I)=0$ and $\mathcal{T}$\ differentiates $L^{1}$.
\end{definition}

By removing the measure zero exceptional set $E(\mathcal{T})=\bigcup
_{I\in\mathcal{T}}\bigcup_{\substack{J_{1},J_{2}\in\mathcal{C}(I)\\J_{1}\neq
J_{2}}}(J_{1}\cap J_{2})$ we may replace the almost disjointness above by disjointness.

Now given any tree $\mathcal{T}$ we define the maximal operator associated to
it as follows
\begin{equation}
M_{\mathcal{T}}\phi(x)=\sup\left\{  \frac{1}{\mu(I)}\int_{I}\left\vert
\phi\right\vert d\mu:x\in I\in\mathcal{T}\right\}  \label{t1}%
\end{equation}
for every $\phi\in L^{1}(X,\mu)$.

The above setting can be used not only for the dyadic maximal operator but
also for the maximal operator on martingales, hence many of the results here
can be viewed as generalizations and refinements of the classical Doob's inequality.

Also for any locally integrable positive function $w$ on $X$, which will be
called weight, we denote $\sigma=w^{-\frac{1}{p-1}}$, and for any
$I\in\mathcal{T}$ we write $w(I)=\int_{I}wd\mu$, $\sigma(I)=\int_{I}\sigma
d\mu$. Now we give the following.

\begin{definition}
A weight $w$ on $X$ will be called $A_{p}$ with respect to $\mathcal{T}$\ if
the following expression%
\[
\lbrack w]_{\mathcal{T}\text{,}p}=[w]_{p}=\sup_{I\in\mathcal{T}}%
\frac{w(I)\sigma(I)^{p-1}}{\mu(I)^{p}}%
\]
is finite.
\end{definition}

A way to study estimates for the above maximal operator is through the
symmetrization of $\phi$ as has been introduced in \cite{Mel2} and \cite{Nik}
and used in \cite{Mel5} to evaluate Bellman functions related to Lorentz
norms. In order to apply this in the context of weights we introduce the
following condition on a weight $w$ on $X$.

\begin{definition}
A weight $w$ on $X$ will be called $A_{p}^{\ast}$ if for some equimeasurable
rearrrangement $w^{\ast\ast}$ of $w$ on $(0,1)$ (not necessarily decreasing)
there exist two constants $c,a>0$ such that for every $t$ in $(0,1]$ the
following estimate holds%
\begin{equation}
\int_{t}^{1}\frac{w^{\ast\ast}(s)}{s^{p}}ds+c\leq a\frac{w^{\ast\ast}%
(t)}{t^{p-1}}\label{Ap1}%
\end{equation}
and also%
\begin{equation}
\lim_{t\rightarrow0^{+}}t^{p}\int_{t}^{1}\frac{w^{\ast\ast}(s)}{s^{p}%
}ds=0\label{Ap1a}%
\end{equation}

\end{definition}

\bigskip Note that by writing $r(t)=\frac{w^{\ast\ast}(t)}{t^{p-1}}$ the first
condition implies that $r(t)>\frac{c}{a}>0$ for all $t$ hence $\lim
_{t\rightarrow0^{+}}\int_{t}^{1}\frac{w^{\ast\ast}(s)}{s^{p}}ds=+\infty$ and
so  $\lim_{t\rightarrow0^{+}}r(t)=0$. Hence we conclude that there is a best
possible pair $(a,c)$ for each such weight, namely $a=\sup_{t}r(t)^{-1}%
\int_{t}^{1}\frac{r(t)}{t}dt$ and $c=\sup_{t}(ar(t)-\int_{t}^{1}\frac{r(t)}%
{t}dt)$. We will refer to this pair as the \textit{constants of the
corresponding }$A_{p}^{\ast}$ weight $w$.

\bigskip\textbf{Example. }Suppose that $w^{\ast\ast}(t)=kt^{b}$ with
$k,b\in\mathbf{R}$. Then the above conditions hold if and only if $-1<b<p-1$
which is exactly the range making $w^{\ast\ast}$ an $A_{p}$ weight on $(0,1)$.
Moreover the corresponding constants $c,a$ can be easily seen to be
$a=\frac{1}{p-1-b},~c=\frac{k}{p-1-b}$.

Now we take into consideration the following theorem proved in \cite{Nik} and
\cite{Mel2}.

\begin{theorem}
Let $G:[0,+\infty)\rightarrow\lbrack0,+\infty)$ be non-decrasing,
$h:(0,1]\rightarrow$ $\mathbb{R}^{+}$ be any locally integrable function. Then
for any nonatomic probability space $(X,\mu)$, equipped with any tree-like
family $\mathcal{T}$ , for any non-increasing right continuous integrable
function $g:(0,1]\rightarrow$ $\mathbb{R}^{+}$ and any $k\in(0,1]$, the
following equality holds (where by $\psi^{\ast}$ we denote the decreasing
equimeasurable rearrangement of $\psi$):en
\begin{gather*}
\sup\left\{  \int_{0}^{k}G[(M_{\mathcal{T}}\phi)^{\ast}(t)]h(t)dt:\phi\text{
measurable on }X\text{ with }\phi^{\ast}=g\right\}  =\\
=\int_{0}^{k}G\left(  \frac{1}{t}\int_{0}^{t}g(u)du\right)  h(t)dt.
\end{gather*}

\end{theorem}

After this given an $A_{p}^{\ast}$ weight $w$, we define the following variant
of the Bellman function (\ref{i6a}).%
\begin{equation}
\mathcal{B}_{p,w}^{\ast}(F,f)=\sup{\huge \{}\int_{0}^{1}((M_{\mathcal{T}}%
\phi)^{\ast})^{p}w^{\ast\ast}:\int_{0}^{1}(\phi^{\ast})^{p}w^{\ast\ast}%
=F,\int_{X}\phi=f{\huge \}}%
\end{equation}
where here by $\phi^{\ast},~(M_{\mathcal{T}}\phi)^{\ast}$ we denote the
equimeasurable \textbf{decreasing }rearrangement of $\phi$ and $M_{\mathcal{T}%
}\phi$ whereas by $w^{\ast\ast}$ we denote the equimeasurable rearrangement of
$w$ that appears in the above definition. Note that in case $w^{\ast\ast}$ is
decreasing $\int_{X}((M_{\mathcal{T}}\phi)^{\ast})^{p}w^{\ast\ast}$ is greater
than or equal to $\int_{X}(M_{\mathcal{T}}\phi)^{p}w$ and $\int_{0}^{1}%
(\phi^{\ast})^{p}w^{\ast\ast}$ is greater than or equal to $\int_{X}\phi^{p}w$
and when $w^{\ast\ast}$ is increasing then the opposite relations hold. Then
we can prove the following.

\begin{theorem}
For the above function we have%
\[
\mathcal{B}_{p,w}^{\ast}(F,f)=(p-1)^{p}a^{p}F\omega_{p}\left(  \frac{cf^{p}%
}{(p-1)^{p-1}a^{p}F}\right)  ^{p}.
\]
where $c,a$ are the constants of the $A_{p}^{\ast}$ weight $w$, the domain of
this function being all $(F,f)$ such that $cf^{p}\leq(p-1)^{p-1}a^{p}F$.

\begin{proof}
In view of the above mentioned result it suffices to consider the expression
$\Delta_{w}(g)=\int_{0}^{1}(t^{-1}\int_{0}^{t}g(u)du)^{p}w^{\ast\ast}(t)dt$
when $g$ runs over all nonnegative decreasing right continuous functions on
$(0,1]$ satisfying $\int_{0}^{1}g(t)dt=f$ and $\int_{0}^{1}g(t)^{p}w^{\ast
\ast}(t)dt=F$. We next define the following function on $(0,1)$%
\[
u(t)=\int_{t}^{1}\frac{w^{\ast\ast}(s)}{s^{p}}ds+c
\]
so that $u^{\prime}(t)=t^{-p}w^{\ast\ast}(t)$. Considering first any bounded
such function $g$ we compute by integration by parts%
\begin{gather*}
\int_{0}^{1}u(t)(\int_{0}^{t}g(u)du)^{p-1}g(t)dt=\frac{1}{p}\int_{0}%
^{1}u(t)[(\int_{0}^{t}g(u)du)^{p}]^{\prime}dt=\\
=\frac{1}{p}(\int_{0}^{1}g(u)du)^{p}u(1)+\frac{1}{p}\int_{0}^{1}(t^{-1}%
\int_{0}^{t}g(u)du)^{p}w^{\ast\ast}(t)=c\frac{f^{q}}{p}+\frac{1}{p}\Delta(g)
\end{gather*}
the integration by parts term $\lim_{t\rightarrow0+}u(t)(\int_{0}%
^{t}g(u)du)^{p}$ being zero because of condition (\ref{Ap1a}) since $g$ is
assumed bounded. Now using Young's inequality $xy\leq\frac{xp}{p}%
+\frac{y^{p^{\prime}}}{p^{\prime}}$ (where $p^{\prime}=p/(p-1)$) in the first
integral as follows, ($\lambda>0$ to be determined later) combined with the
condition $\frac{u(t)t^{p-1}}{w^{\ast\ast}(t)}\leq a$ from the above
definition we get%
\begin{gather*}
\int_{0}^{1}u(t)(\int_{0}^{t}g(u)du)^{p-1}g(t)dt=\\
=\int_{0}^{1}(\lambda g(t)w^{\ast\ast}(t)^{1/p})(\frac{w^{\ast\ast}(t)^{1/p}%
}{\lambda^{1/(p-1)}t}\int_{0}^{t}g(u)du)^{p-1}\frac{u(t)t^{p-1}}{w^{\ast\ast
}(t)}dt\leq\\
\leq\frac{a}{p}\int_{0}^{1}\lambda^{p}g(t)^{p}w^{\ast\ast}(t)dt+\frac
{a}{p^{\prime}}\int_{0}^{1}\lambda^{-p^{\prime}}(\frac{1}{t}\int_{0}%
^{t}g(u)du)^{p}w^{\ast\ast}(t)dt=\\
=\frac{a\lambda^{p}}{p}\int_{0}^{1}g(t)^{p}w^{\ast\ast}(t)dt+a\frac
{\lambda^{-p^{\prime}}}{p^{\prime}}\int_{0}^{1}(\frac{1}{t}\int_{0}%
^{t}g(u)du)^{p}w^{\ast\ast}(t)dt=a\frac{\lambda^{p}}{p}F+a\frac{\lambda
^{-p^{\prime}}}{p^{\prime}}\Delta_{w}(g)\text{.}%
\end{gather*}
Therefore we have by writing $\lambda^{p^{\prime}}=(p-1)a(\beta+1),~\beta>0$
and using the above inequalities we get that%
\begin{equation}
\Delta_{w}(g)\leq(1+\frac{1}{\beta})\frac{(\beta+1)^{p-1}(p-1)^{p}%
a^{p}F-(p-1)cf^{p}}{(p-1)}.\label{Ap2}%
\end{equation}
Next, given an arbitrary $g$, the above estimate can be used for the
truncations $g_{M}=\min(g,M)$ and $F,f$ replaced by the corresponding
quantities for $g_{M}$ and then take $M\rightarrow+\infty$ and use monotone
convergence to infer that (\ref{Ap2}) holds for the general nonnegative
decreasing right continuous function on $(0,1]$ satisfying $\int_{0}%
^{1}g(t)dt=f$ and $\int_{0}^{1}g(t)^{p}w^{\ast\ast}(t)dt=F$. Moreover since
$\Delta_{w}(g)>0$ the inequality  (\ref{Ap2}) implies that $(\beta
+1)^{p-1}(p-1)^{p}a^{p}F-(p-1)cf^{p}>0$ for every $\beta>0$ and so letting
$\beta\rightarrow0^{+}$ we conclude that $(F,f)$ must satisfy the inequality
$cf^{p}\leq(p-1)^{p-1}a^{p}F$ given in the statement of the Theorem.

Writing $A=(p-1)^{p}a^{p}F$ and $B=(p-1)cf^{p}$ it is easy to compute (see for
example \cite{Mel1} pg. 326) that the minimum possible value of the right hand
side of (\ref{Ap2}) is equal to $A\omega_{p}\left(  \frac{B}{A}\right)  ^{p}$.
This proves the inequality
\begin{equation}
\mathcal{B}_{p,w}^{\ast}(F,f)\leq(p-1)^{p}a^{p}F\omega_{p}\left(  \frac
{cf^{p}}{(p-1)^{p-1}a^{p}F}\right)  ^{p}.\label{Ap2a}%
\end{equation}
Now we consider the continuous positive decreasing function%
\begin{equation}
g_{\alpha}(t)=f(1-\alpha)t^{-\alpha}\label{Ap3}%
\end{equation}
where $0\leq\alpha<1$, and any $A_{p}^{\ast}$ weight $w$ that is
equimeasurable to%
\begin{equation}
w^{\ast\ast}(t)=kt^{b}\text{, }k>0,-1<b<p-1\label{Ap4}%
\end{equation}
Clearly $\int_{0}^{1}g_{\alpha}(t)dt=f$ and $\int_{0}^{1}g_{\alpha}%
(t)^{p}w^{\ast\ast}(t)dt=\frac{kf^{p}(1-\alpha)^{p}}{1+b-\alpha p}$ assuming
that $\alpha<\frac{1+b}{p}$. Next note that $\frac{1}{t}\int_{0}^{t}g_{\alpha
}(u)du=\frac{g_{\alpha}(t)}{1-\alpha}$ for all $t\in(0,1]$ and so we have
$\Delta_{w}(g_{\alpha})=\left(  \frac{1}{1-\alpha}\right)  ^{p}\int_{0}%
^{1}g_{\alpha}(t)^{p}w^{\ast\ast}(t)dt$. The condition $\int_{0}^{1}g_{\alpha
}(t)^{p}w^{\ast\ast}(t)dt=F$ is then equivalent to the following equation in
$\alpha$%
\begin{equation}
\frac{(1-\alpha)^{p}}{1+b-\alpha p}=\frac{F}{kf^{q}}\text{.}\label{Ap5}%
\end{equation}
To study this equation we write%
\begin{equation}
z=\frac{p-1-b}{p-1}\frac{1}{1-\alpha}\label{Ap6}%
\end{equation}
and note that (\ref{Ap5}) is then equivalent to%
\begin{equation}
-(p-1)z^{p}+pz^{p-1}=\frac{kf^{p}}{(\frac{p-1}{p-1-b})^{p-1}F}%
\end{equation}
thus%
\begin{equation}
z=\omega_{p}(\frac{kf^{p}}{(\frac{p-1}{p-1-b})^{p-1}F})\label{Ap7}%
\end{equation}
and so using (\ref{Ap6})%
\begin{equation}
\Delta_{w}(g_{\alpha})=(\frac{p-1}{p-1-b})^{p}F\omega_{p}(\frac{kf^{p}}%
{(\frac{p-1}{p-1-b})^{p-1}F}).
\end{equation}
But now note that the constants $c,a$ of the weight $w$ are $a=\frac{1}%
{p-1-b},~c=\frac{k}{p-1-b}$ and so%
\begin{equation}
\Delta_{w}(g_{\alpha})=(p-1)^{p}a^{p}F\omega_{p}\left(  \frac{cf^{p}%
}{(p-1)^{p-1}a^{p}F}\right)  ^{p}%
\end{equation}
and moreover by varying $k,b$ with $-1<b<p-1$ we can achieve all possible
pairs of constants $c,a$. {}This completes the proof.
\end{proof}
\end{theorem}

\section{Estimation via double maximal operators}

Here we will use an inequality introduced by A. Lerner, see \cite{Ler}, for
the nondyadic case. We fix $p>1$, let $w$ be an $A_{p}$ weight with respect to
the tree $\mathcal{T}$ and we denote for any $I$ in $\mathcal{T}$,
$w(I)=\int_{I}wd\mu$, $\sigma=w^{-\frac{1}{p-1}}$, $\sigma(I)=\int_{I}\sigma
d\mu$. Also by $M_{\mathcal{T},w}$ we denote the maximal operator with respect
to the tree $\mathcal{T}$ \ but when $X$ is equipped by the measure $w\mu$
instead of $\mu$, and similarly for $M_{\mathcal{T}\text{,}\sigma}$. Then the
following holds.

\begin{proposition}
\bigskip Let $w$ be an $A_{p}$ weight with respect to the tree $\mathcal{T}$
and $\mathcal{T}$-constant $[w]_{p}=\sup_{I\in\mathcal{T}}\frac{w(I)\sigma
(I)^{p-1}}{\mu(I)^{p}}$ Then for any $\phi$ we have the following pointwise
estimate%
\[
(M_{\mathcal{T}}\phi)^{p-1}\leq\lbrack w]_{p}M_{\mathcal{T},w}[(M_{\mathcal{T}%
\text{,}\sigma}(\phi\sigma^{-1}))^{p-1}w^{-1}].
\]

\begin{proof}
The proof follows from the following inequalities valid for any $I\in
\mathcal{T}$.%
\begin{align*}
\left(  \frac{1}{\mu(I)}\int_{I}\phi d\mu\right)  ^{p-1}  &  =\frac
{w(I)\sigma(I)^{p-1}}{\mu(I)^{p}}\left(  \frac{\mu(I)}{w(I)}\left(  \frac
{1}{\sigma(I)}\int_{I}\phi\sigma^{-1}\sigma d\mu\right)  ^{p-1}\right)  \leq\\
&  \leq[w]_{p}\frac{1}{w(I)}\int_{I}M_{\mathcal{T}\text{,}\sigma}(\phi
\sigma^{-1})^{p-1}w^{-1}wd\mu
\end{align*}
since $M_{\mathcal{T}\text{,}\sigma}(\phi\sigma^{-1})(x)\geq\frac{1}%
{\sigma(I)}\int_{I}\phi\sigma^{-1}\sigma d\mu$ for every $x$ in $I$.
\end{proof}
\end{proposition}

\bigskip As a first application of this fixing a tree $\mathcal{T}$ on a
probability space $(X,\mu)$ and given an $A_{p}$ weight $w$ in the sense of
Definition 2, we define the following generalization of the Bellman function
(\ref{i6a}), where $p>1$%
\begin{equation}
\mathcal{B}_{p,w}^{\mathcal{T}}(F,f)=\sup{\huge \{}\int_{X}(M_{\mathcal{T}%
}\phi)^{p}wd\mu:\int_{X}\phi^{p}wd\mu=F,\int_{X}\phi d\mu=f{\huge \}}%
\end{equation}
and we have the following estimates

\begin{theorem}
For any tree $\mathcal{T}$ on a probability space $(X,\mu)$ and any $A_{p}$
weight $w$ and any $\phi$ with $\int_{X}\phi^{p}wd\mu=F,\int_{X}\phi d\mu=f$
we have%
\begin{gather}
\int_{X}(M_{\mathcal{T}}\phi)^{p}wd\mu\leq\nonumber\\
\leq\lbrack w]_{p}^{1/(p-1)}F\omega_{p}\left(  \frac{f^{p}}{\sigma(X)^{p-1}%
F}\right)  ^{p}\omega_{p^{\prime}}\left(  \frac{(\int_{X}\phi^{p-1}%
wd\mu)^{p^{\prime}}}{w(X)^{p^{\prime}-1}F\omega_{p}\left(  \frac{f^{p}}%
{\sigma(X)^{p-1}F}\right)  ^{p}}\right)  ^{p^{\prime}}\label{W1}%
\end{gather}
In particular%
\begin{equation}
\mathcal{B}_{p,w}^{\mathcal{T}}(F,f)\leq p^{p^{\prime}}[w]_{p}^{1/(p-1)}%
F\omega_{p}\left(  \frac{f^{p}}{\sigma(X)^{p-1}F}\right)  ^{p}\label{W2}%
\end{equation}

\begin{proof}
By applying estimate (1.12) after Theorem 1 in \cite{Mel1} for the exponent
$p^{\prime}=\frac{p}{p-1}$ to the function $\rho=(M_{\mathcal{T}\text{,}%
\sigma}(\phi\sigma^{-1}))^{p-1}w^{-1}$ and with respect to the tree
$\mathcal{T}$ but on the probability space $(X,\frac{1}{w(X)}wd\mu)$ (where as
usual $w(X)=\int_{X}wd\mu$) we get%
\begin{gather}
\frac{1}{[w]_{p}^{1/(p-1)}}\int_{X}(M_{\mathcal{T}}\phi)^{p}wd\mu\leq
w(X)\int_{X}(M_{\mathcal{T},w}\rho)^{p^{\prime}}w\frac{d\mu}{w(X)}%
\leq\nonumber\\
\leq w(X)\int_{X}\rho^{p^{\prime}}w\tfrac{d\mu}{w(X)}.\omega_{p^{\prime}%
}\left(  \frac{(\int_{X}\rho w\tfrac{d\mu}{w(X)})^{p^{\prime}}}{\int_{X}%
\rho^{p^{\prime}}w\tfrac{d\mu}{w(X)}}\right)  ^{p^{\prime}}.
\end{gather}
Note that (as proved in  \cite{Mel1}) the function $x\omega_{p^{\prime}%
}\left(  \dfrac{y^{p^{\prime}}}{x}\right)  ^{p^{\prime}}$ is increasing in $x$
and decreasing in $y$. Now we have%
\[
\int_{X}\rho w\tfrac{d\mu}{w(X)}=\int_{X}(M_{\mathcal{T}\text{,}\sigma}%
(\phi\sigma^{-1}))^{p-1}\tfrac{d\mu}{w(X)}\geq\int_{X}(\phi\sigma^{-1}%
)^{p-1}\tfrac{d\mu}{w(X)}=\int_{X}\phi^{p-1}w\tfrac{d\mu}{w(X)}%
\]
and using estimate (1.12) after Theorem 1 in \cite{Mel1} for the exponent $p$
to the function $\rho=\phi\sigma^{-1}$ and with respect to the tree
$\mathcal{T}$ but on the probability space $(X,\frac{1}{\sigma(X)}\sigma
d\mu)$ we get (since $\sigma^{-(p-1)}=w$)%
\begin{gather}
\int_{X}\rho^{p^{\prime}}wd\mu=\int_{X}(M_{\mathcal{T}\text{,}\sigma}%
(\phi\sigma^{-1}))^{p}w^{-p^{\prime}}.wd\mu=\int_{X}(M_{\mathcal{T}%
\text{,}\sigma}(\phi\sigma^{-1}))^{p}\sigma d\mu\leq\nonumber\\
\sigma(X)\int_{X}(\phi\sigma^{-1})^{p}\sigma\tfrac{d\mu}{\sigma(X)}.\omega
_{p}\left(  \frac{(\int_{X}\phi\sigma^{-1}\sigma\tfrac{d\mu}{\sigma(X)})^{p}%
}{\int_{X}(\phi\sigma^{-1})^{p}\sigma\tfrac{d\mu}{\sigma(X)}}\right)
^{p}=\nonumber\\
=\int_{X}\phi^{p}wd\mu.\omega_{p}\left(  \frac{(\int_{X}\phi d\mu)^{p}}%
{\sigma(X)^{p-1}\int_{X}\phi^{p}wd\mu}\right)  ^{p}=F\omega_{p}\left(
\frac{f^{p}}{\sigma(X)^{p-1}F}\right)  ^{p}\text{.}%
\end{gather}
Now combining the above estimates we get%
\begin{gather}
\frac{1}{[w]_{p}^{1/(p-1)}}\int_{X}(M_{\mathcal{T}}\phi)^{p}wd\mu
\leq\nonumber\\
\leq F\omega_{p}\left(  \frac{f^{p}}{\sigma(X)^{p-1}F}\right)  ^{p}%
\omega_{p^{\prime}}\left(  \frac{(\int_{X}\phi^{p-1}wd\mu)^{p^{\prime}}%
}{w(X)^{p^{\prime}-1}F\omega_{p}\left(  \frac{f^{p}}{\sigma(X)^{p-1}F}\right)
^{p}}\right)  ^{p^{\prime}}%
\end{gather}
which proves (\ref{W1}). Since $\omega_{p^{\prime}}(x)\leq\frac{p^{\prime}%
}{p^{\prime}-1}=p$ the estimate (\ref{W2}) follows also.
\end{proof}
\end{theorem}

\bigskip To get lower bounds for the Bellman function we invoke the following construction.

Fixing $\alpha$ with $0<\alpha<1$ and using Lemma 1 in  \cite{Mel1}, we fix
now a tree $\mathcal{T}$, for example the dyadic subintervals of $[0,1]$, and
choose for every $I\in\mathcal{T}$ a family $\mathcal{F}(I)\subseteq
\mathcal{T}$ of pairwise almost disjoint subsets of $I$ such that
\begin{equation}
\sum_{J\in\mathcal{F}(I)}\mu(J)=(1-\alpha)\mu(I)\text{.}\label{e13}%
\end{equation}
Then we define $\mathcal{S}=\mathcal{S}_{\alpha}$ to be the smallest subset of
$\mathcal{T}$ such that $X\in\mathcal{S}$ and for every $I\in\mathcal{S}$,
$\mathcal{F}(I)\subseteq\mathcal{S}$. Next for every $I\in\mathcal{S}$ we
define the set
\begin{equation}
A_{I}=I~\backslash%
%TCIMACRO{\dbigcup \limits_{J\in\mathcal{F}(I)}}%
%BeginExpansion
{\displaystyle\bigcup\limits_{J\in\mathcal{F}(I)}}
%EndExpansion
J\label{e14}%
\end{equation}
and note that $\mu(A_{I})=\alpha\mu(I)$ and $I=%
%TCIMACRO{\dbigcup \limits_{_{\substack{J\in\mathcal{S}\\J\subseteq I}}}}%
%BeginExpansion
{\displaystyle\bigcup\limits_{_{\substack{J\in\mathcal{S}\\J\subseteq I}}}}
%EndExpansion
A_{J}$ for every $I\in\mathcal{S}$. Also since $\mathcal{S}=\bigcup_{m\geq
0}\mathcal{S}_{(m)}$ where $\mathcal{S}_{(0)}=\{X\}$ and $\mathcal{S}%
_{(m+1)}=\bigcup_{I\in\mathcal{S}_{(m)}}\mathcal{F}(I)$, we can define
rank$(I)=r(I)$ for $I\in\mathcal{S}$ to be the unique integer $m$ such that
$I\in\mathcal{S}_{(m)}$ and remark that $\sum\limits_{\substack{\mathcal{S}\ni
J\subseteq I\\r(J)=r(I)+m}}\mu(J)=(1-\alpha)^{m}\mu(I)$ for every
$I\in\mathcal{S}$.

Next for any $\lambda,\gamma>0$ we define the function%
\begin{equation}
\psi=%
%TCIMACRO{\dsum \limits_{I\in\mathcal{S}}}%
%BeginExpansion
{\displaystyle\sum\limits_{I\in\mathcal{S}}}
%EndExpansion
\lambda\gamma^{r(I)}\chi_{A_{I}}%
\end{equation}
and we have for any $I\in\mathcal{S}$ the following%
\begin{equation}
\frac{1}{\mu(I)}\int_{I}\psi d\mu=\frac{\lambda\alpha}{1-\gamma(1-\alpha
)}\gamma^{r(I)}\text{.}%
\end{equation}
Hence taking
\begin{equation}
\phi_{\alpha}=%
%TCIMACRO{\dsum \limits_{I\in\mathcal{S}}}%
%BeginExpansion
{\displaystyle\sum\limits_{I\in\mathcal{S}}}
%EndExpansion
\lambda_{1}\gamma_{1}^{r(I)}\chi_{A_{I}}\text{, \ }w_{\alpha}=%
%TCIMACRO{\dsum \limits_{I\in\mathcal{S}}}%
%BeginExpansion
{\displaystyle\sum\limits_{I\in\mathcal{S}}}
%EndExpansion
\lambda_{2}\gamma_{2}^{r(I)}\chi_{A_{I}}%
\end{equation}
we have for any $I\in\mathcal{S}$
\begin{equation}
\frac{w_{\alpha}(I)[w_{\alpha}^{-\frac{1}{p-1}}(I)]^{p-1}}{\mu(I)^{p}}%
=\frac{\alpha^{p}}{[1-\gamma_{2}(1-\alpha)][1-\gamma_{2}^{-\frac{1}{p-1}%
}(1-\alpha)]}%
\end{equation}
thus $w_{\alpha}$ is an $A_{p}$ weight but with respect to the \textit{tree
}$\mathcal{S}_{\alpha}$ on $(X,\mu)$ and with $[w_{\alpha}]_{p}$ equal to the
right hand side of the above relation. Moreover%
\begin{equation}
M_{\mathcal{S}}\phi_{\alpha}\geq%
%TCIMACRO{\dsum \limits_{I\in\mathcal{S}}}%
%BeginExpansion
{\displaystyle\sum\limits_{I\in\mathcal{S}}}
%EndExpansion
\frac{1}{\mu(I)}\int_{I}\phi_{a}d\mu\chi_{A_{I}}=\frac{\alpha}{1-\gamma
(1-\alpha)}\phi_{\alpha}\text{.}%
\end{equation}
However the values of such functions on each $A_{I}$ where $r(I)=m$ is of the
form
\begin{equation}
\gamma_{m}=\frac{\lambda}{\alpha(1-\alpha)^{m}}\int_{(1-\alpha)^{m+1}%
}^{(1-\alpha)^{m}}u^{s}du\label{e16}%
\end{equation}
for some real numbers $\lambda,s>0$ and as it is proved in Lemma 3 of
\cite{Mel5} these behave like functions of the form $\lambda t^{s}$ on $(0,1]$
as we approach the limit $\alpha\rightarrow0^{+}$. Hence by taking a sequence
$\alpha_{m}\rightarrow0$ considering the trees $\mathcal{T}_{m}=\mathcal{S}%
_{\alpha_{m}}$ on $(X,\mu)$ and using the construction for the lower bound in
the proof of Theorem 2, choosing the constants $k,b$ ($-1<b<p-1$)
appropriately in  (\ref{Ap4}) according to the conditions $a=\frac{1}{p-1-b}$,
$c=ka$, $\frac{k}{b+1}=z$, $\frac{1}{b+1}(\frac{p-1}{p-1-b})^{p-1}=h$ from the
restrictions below which give $\frac{p-1-b}{p-1}=\omega_{p}(\frac{1}{h})$ we
conclude the following.

\begin{proposition}
Given appropriate $F,f,h,z$ there exists a sequence of trees $\mathcal{T}_{m}$
on $(X,\mu)~$and two sequences $(\phi_{m})$ and $(w_{m})$ of positive
measurable functions on $(X,\mu)$ such that $\int_{X}\phi_{m}d\mu\rightarrow
f$,~ $\int_{X}\phi_{m}^{p}w_{m}d\mu\rightarrow F$, each $w_{m}$ is an $A_{p}$
weight with respect to the tree $\mathcal{T}_{m}$ with $[w_{m}]_{p}\rightarrow
h$ and $\int_{X}w_{m}d\mu\rightarrow z$ such that%
\begin{equation}
\lim_{m\rightarrow\infty}\int_{X}(M_{\mathcal{T}_{m}}\phi_{m})^{p}w_{m}%
d\mu\geq F\omega_{p}\left(  \frac{zf^{p}}{hF}\right)  ^{p}\omega_{p}(\frac
{1}{h})^{-p}\text{.}%
\end{equation}

\end{proposition}

The above proposition implies a lower bound on the class of functions
$\mathcal{B}_{p,w}^{\mathcal{T}}(F,f)$ when viewed over all trees
$\mathcal{T}$ and corresponding $A_{p}$ weights $w$.

\end{document}